
\documentstyle{amsppt}
\mag=\magstep0
\hsize=6.5 true in
\vsize=9.0 true in
\voffset=-0.5 true in

\define\normalskip{\baselineskip=18pt}
\define\tableskip{\baselineskip=12pt}
\define\Hecke{{\Cal H}}
\define\complexes{{\Bbb C}}
\define\Irr{\operatorname{Irr}}
\define\Ind{\operatorname{Ind}}
\define\Poincare{{Poincar\'e}}
\define\qbinomial#1#2{\text{\bf C}^{#1}_{#2}}
\define\abs#1{\left\vert{#1}\right\vert}
\define\trace{\operatorname{trace}}
\define\residue#1#2{{\left [ #1 \right ]_{#2}}}
\define\setof#1#2{\left\{#1\,\vert\,#2\right\}}


\define\TypeATheorem{Theorem~1.1}
\define\TypeACorollary{Corollary~1.2}
\define\DihedralTheorem{Theorem~1.3}
\define\CoxeterTheorem{Theorem~1.4}

\define\PalindromicProposition{Proposition~2.1}
\define\DegreeCorollary{Corollary~2.2}
\define\RecursionProposition{Proposition~2.3}
\define\ProductProposition{Proposition~2.4}

\define\DihedralLemma{Lemma~4.1}

\define\ExampleB{Example~6.1}
\define\ExampleD{Example~6.2}
\define\ExampleF{Example~6.3}
\define\ExampleH{Example~6.4}
\define\Conjecture{Conjecture~6.5}


\define\AlvisLusztigRef{1}
\define\AlvisRatioRef{2}
\define\BensonRef{3}
\define\BourbakiRef{4}
\define\CarterClassRef{5}
\define\CarterSimpleGroupsRef{6}
\define\CarterFiniteGroupsRef{7}
\define\FrameOneRef{8}
\define\FrameTwoRef{9}
\define\GeckPfeifferRef{10}
\define\GreenRef{11}
\define\GroveRef{12}
\define\GyojaRef{13}
\define\KacCheungRef{14}
\define\KilmoyerSolomonRef{15}
\define\KondoRef{16}
\define\LusztigClassicalRef{17}
\define\LusztigUnipRef{18}
\define\LusztigBensonCurtisRef{19}
\define\LusztigClassTwoRef{20}
\define\SchonertRef{21}
\define\StanleyRef{22}
\define\SurowskiOneRef{23}
\define\SurowskiTwoRef{24}


\document

\headline={\hss \tenrm --\  \folio\ -- \hss}
\footline={\hfill}


\head 
\centerline
{\bf A family of polynomials related to generic degrees} \endhead

\bigskip
\centerline{Dean Alvis}
\centerline{Department of Mathematical Sciences}
\centerline{Indiana~University~South~Bend}
\centerline{South~Bend, IN, 46634, USA}
\centerline{e-mail: \ \ dalvis\@iusb\.edu}
\bigskip

\centerline{
\vbox{\hsize 5 true in \noindent
\eightrm
Abstract. \ 
A family of polynomials parameterized by the 
conjugacy classes of a finite Coxeter group is
investigated.  
These polynomials, together with the
character table of the group, 
determine the associated 
 generic degrees.   
The polynomials are described  
completely for classes that meet a 
parabolic subgroup whose components
are of type A or are dihedral, and 
for the class of Coxeter elements.
}}

\normalskip


\head
\leftline{\bf 1. Introduction}
\endhead

Let $W$ be a finite Coxeter
group with set of distinguished generators $S$,
length function $\ell:w \mapsto \ell(w)$, and
set of irreducible complex characters $\Irr(W)$.
Let $K$ be a splitting field for $W$ contained in the
field of complex numbers, 
let $q$ be an indeterminate over $K$,
and let $\Hecke$ be the
Hecke algebra of $W$ over $K(\sqrt{q})$ with 
standard basis elements $T_w$, $w \in W$.
Thus
$$
\cases
T_s^2 = q T_e + (q-1)T_s, &
\text{if $s \in S$,} \cr
T_x T_y = T_{xy}, &
\text{if $\ell(x)+\ell(y)=\ell(xy)$},\cr
\endcases
$$
where $e$ is the identity of $W$.
It is known that $\Hecke$
is split semisimple over $K(\sqrt{q})$, and that
$\widetilde{\chi}(T_w) \in K[\sqrt{q}]$ for
any irreducible character $\widetilde{\chi}$ of
$\Hecke$ and $w \in W$.
Moreover, if $\chi(w)$ is obtained from
$\widetilde{\chi}(T_w)$ by the substitution
$\sqrt{q}\mapsto 1$
for $w \in W$, then
the correspondence
$\chi \leftrightarrow \widetilde{\chi}$
is a bijection between $\Irr(W)$ 
and the set of
irreducible characters of $\Hecke$
(\cite {\GyojaRef},
\cite{\LusztigBensonCurtisRef},
\cite{\AlvisLusztigRef}).
The {\it generic degrees}
$\left\{d_\chi(q): \chi \in \Irr(W)\right\}$
are determined
by
$$
\sum
\limits_{\chi \in \Irr(W)}
\widetilde{\chi}(T_w)
d_\chi(q) =
\cases
P(q), & \text{if $w = e$,} \cr
0, & \text{if $w \ne e$,} \cr
\endcases
\leqno{(1)}
$$
where 
$P(q) = P_W(q)= \sum_{w \in W} q^{\ell(w)}$
is the {\it {\Poincare} polynomial} of $W$
(\cite{\AlvisLusztigRef}).  
The generic degrees are important 
in the representation theory of
finite groups of Lie type because they 
appear in the formulas for the
degrees of the irreducible
representations of such groups.

This paper is an investigation
of the polynomials
$$
g_w (q) = \sum_{\chi \in \Irr(W)} \chi(w) d_\chi(q),
\qquad\quad
w \in W,
\leqno{(2)}
$$
obtained by applying the substitution
$\sqrt{q}\mapsto 1$ to the first factor
in each summand of (1).
Clearly $g_w(q) = g_{w^\prime}(q)$ if 
$w$, $w^\prime$ are conjugate in $W$.  
If $\sigma$ is an automorphism of $K$, then
$d_{{}^\sigma\chi}(q)$
can be obtained by applying $\sigma$ to the
coefficients of $d_\chi(q)$,
and it follows that 
$g_w(q)$ has rational coefficients
for all $w \in W$.
The orthogonality relations for $W$
combine with (2) to give
$$
d_\chi(q) =
{1 \over \abs{W}}
\sum_{w \in W}
\chi(w^{-1}) g_w(q).
$$
Thus  
the generic degrees are determined by
the polynomials 
$\left\{g_w(q):w \in W\right\}$ together
with the character
table of $W$.


The following results describe $g_w(q)$
for certain $w$ in terms of data that 
can easily be extracted from $w$ and $(W,S)$.

\proclaim{\TypeATheorem}
Suppose $w$ is conjugate to an element of a
standard 
parabolic subgroup of $W$ all of whose
components have type A.  
Let $\ell = \vert S \vert$ be the
rank of $(W,S)$.  
Then
$$
g_w(q) = {(1-q)^\ell P(q) \over \det(1 - q w)},
$$
where $\det(1 - q w)$ is the characteristic polynomial
of $w$ in the natural reflection representation of $W$.
\endproclaim


\proclaim{\TypeACorollary}
We have
$g_e(q)=P(q)$,
and
$g_s(q) = (1-q) P(q) / (1+q)$
for $s \in S$.
\endproclaim


The next result, together with 
{\TypeACorollary}, completely determines
the polynomials
$g_w(q)$ when
$W$ is a 
dihedral group. 

\proclaim{\DihedralTheorem}
Suppose $W$ is dihedral of order $2n$ with
$n \ge 3$, $S=\left\{ r,s \right\}$, and
$w = (rs)^k$ with $1 \le k \le n/2$.
Then
$g_w(q) = (1-q^k)(1-q^{n-k})$.
\endproclaim


Recall that
a {\it Coxeter element} of $W$ is an element
conjugate to a product of the
form $s_1 \cdots s_\ell$,
where
$S = \left \{ s_1, \dots , s_\ell \right \}$.

\proclaim{\CoxeterTheorem}
Let $c$ be a Coxeter element of $W$, and let 
$d_1$, \dots, $d_\ell$ be the degrees of the
basic polynomial invariants of $W$.
Then
$$
g_c(q) = (1-q^{d_1-1}) (1-q^{d_2-1}) 
\cdots (1-q^{d_\ell-1}).
$$
\endproclaim

The proof of {\CoxeterTheorem} 
that appears below is
a case-by-case argument based on the classification 
of indecomposable finite Coxeter groups.  
It is the author's hope that 
results similar to those above,
or at least some description of $g_w(q)$ that is 
independent of the generic degrees,
can be obtained
for arbitrary $w \in W$.
It should be noted that 
$g_w(q)$ is not in general equal
to a product
of cyclotomic polynomials, as
can be seen from the examples in the last section.


\head
\leftline{\bf 2. Preliminary Results}
\endhead

Let $N$ be the number of positive roots in a 
root system corresponding to $(W,S)$.
Denote by $\varepsilon$ the sign character
of $W$, 
and put 
$\varepsilon_w
= \varepsilon(w) = (-1)^{\ell(w)}$
for $w \in W$.
The next result shows that the polynomial
$g_w(q)$ is
palindromic if $w$ is even
($\varepsilon_w = 1$) and
skew palindromic if 
$w$ is odd ($\varepsilon_w = -1$).

\proclaim{\PalindromicProposition}
For any $w \in W$, we have 
$q^{N} g_w(1/q) = \varepsilon_w g_w(q)$.
\endproclaim

\demo{Proof}
We recall that for any irreducible 
character $\chi$ of $W$,
$$
d_{\varepsilon \chi}(q) =
q^N d_{\chi}(1/q)
$$
\cite{\GreenRef, Lemma~8.3}.  
Since
$\chi \mapsto \varepsilon \chi$
is a permutation on
$\Irr(W)$,
$$
q^{N} g_w (1/q) 
=
\sum_{\chi \in \Irr(W)}
\chi(w) q^{N} d_\chi(1/q) 
=
\sum_{\chi\in \Irr(W)}
\chi(w) d_{\varepsilon \chi}(q) 
= 
\varepsilon_w g_w(q),
$$
as required.
\qed\enddemo

Since $q$ divides $d_\chi(q)$ for any $\chi \ne 1$
and $d_1(q)=1$, we have
$g_w(0)=1$ for all $w \in W$.  Hence the
next result follows immediately from the proposition.

\proclaim{\DegreeCorollary}
For any $w \in W$, $g_w(q)$ is a polynomial
of degree $N$.   Moreover,
the constant term of 
$g_w(q)$ is $1$ and the term of highest degree
is $\varepsilon_w q^N$.
\endproclaim


For $J \subseteq S$, let $W_J$ be the standard 
parabolic subgroup of $W$ generated by $S$,
with corresponding {\Poincare} polynomial
$P_{W_J}(q)$.  For $x \in W_J$, define 
$$
g_{W_J,x}(q) = \sum_{\varphi\in \Irr(W_J)}
     \varphi(x) d_\varphi(q).
$$
The next result gives a recursion formula
for $g_w(q)$ if $w$ is conjugate to an element
of a standard parabolic subgroup.

\proclaim{\RecursionProposition}
Suppose $J \subseteq S$ and 
$w$ is conjugate in $W$ to 
$x \in W_J$.
Then
$$
g_w(q) = {P_{W}(q) \over P_{W_J}(q)} \ g_{W_J,x}(q).
$$
\endproclaim

\demo{Proof}
Let $\varphi$ be an irreducible character of $W_J$.
For $\chi$ a character of $W$, let
$\left < \chi \vert_{W_J}, \varphi \right>$ denote the
multiplicity of $\varphi$ in the restriction of $\chi$ to $W_J$.
Then it is known that 
$$
\sum_{\chi \in \Irr(W)} 
\left< \chi \vert_{W_J}, \varphi \right> d_{\chi}(q)
=
{P_{W}(q) \over P_{W_J}(q)} d_\varphi(q)
$$
(\cite{\AlvisRatioRef}).  
Since $w$ is conjugate to $x$ in $W$, 
we have 
$g_w(q) = g_x(q)$, 
and hence
$$
\eqalign{
g_w(q) 
= & \ 
\sum_{\chi \in \Irr(W)} \chi(x) d_\chi(q) 
= 
\sum_{\chi \in \Irr(W)} \left( \sum_{\varphi \in \Irr(W_J)} 
   \left< \chi \vert_{W_J}, \varphi \right> 
   \varphi(x) \right) d_\chi(q) \cr
= & \ 
\sum_{\varphi \in \Irr(W_J)}
\varphi(x)
\left( \sum_{\chi \in \Irr(W)} 
\left< \chi \vert_{W_J}, \varphi \right> 
    d_\chi(q) \right ) \cr
= & \ 
\sum_{\varphi \in \Irr(W_J)} \varphi(x) 
{P_{W}(q) \over P_{W_J}(q)}
\, d_\varphi(q) 
= 
{P_{W}(q) \over P_{W_J}(q)}
\, g_{W_J,x}(q), \cr
}
$$
as required.
\qed\enddemo


Suppose $W$ is the product as a Coxeter group
of the parabolic subgroups $W_{J_1}$,
\dots, $W_{J_k}$; that is,
$S$ is the disjoint
union of the pairwise commuting subsets 
$J_1$, \dots, $J_k$.   If $\chi \in \Irr(W)$, 
then $\chi$
has a unique expression as an outer tensor product
$\chi_1 \times \cdots \times \chi_k$
where $\chi_j \in \Irr(W_{J_j})$ 
for $1 \le j \le k$, and moreover, 
$d_\chi(q) = d_{\chi_1}(q) 
\cdot \cdots \cdot d_{\chi_k}(q)$.
Thus the next result follows directly 
from the definition of $g_w(q)$.  

\proclaim{\ProductProposition}
Suppose $S$ is the disjoint
union of the pairwise commuting subsets 
$J_1$, \dots, $J_k$.   Let
$w = w_1 \cdots w_k \in W$, where
$w_i \in W_{J_i}$.   Then 
$$
g_w(q) = g_{W_{J_1},w_1}(q) \cdot \cdots \cdot
                g_{W_{J_k},w_k}(q).
$$
\endproclaim


\head
\leftline{\bf 3. The proof of \TypeATheorem}
\endhead

Replacing $w$ by a conjugate if necessary, we can
assume that $w$ is an element of a standard parabolic
subgroup $W_J$ with components of type A.
Using {\RecursionProposition} and
{\ProductProposition}, we can reduce to the case
in which $W=W_J$ is indecomposable of
type $A$.  Then
$d_\chi(q)=f_\chi(q)$ for
$\chi \in \Irr(W)$,
where
$$
f_\chi(q) =
{ (1-q)^\ell P(q) \over \abs{W}}
\sum_{x \in W}
{\chi(x) \over \det(1 - q x)}
$$
is the associated 
``fake degree'' polynomial
(\cite{\LusztigClassicalRef,2.2}).
Thus 
$$
\eqalign{
g_w(q) 
= & \ 
\sum_{\chi \in \Irr(W)}\chi(w) f_\chi(q) \cr
= & \ 
\sum_{\chi \in \Irr(W)} \chi(w) 
{(1-q)^\ell P(q) \over \abs{W}}
\sum_{x \in W}
{\chi(x) \over \det(1 - q x)}\cr
= & \ 
{(1-q)^\ell P(q) \over \abs{W}}
\sum_{x \in W}
{1 \over \det(1 - q x )}
\left( \sum_{\chi \in \Irr(W)} \chi(w)\chi(x) \right) . \cr
}
$$
Since the characters of $W$ are real-valued, the
orthogonality relations for $W$ show that the 
inner sum vanishes unless $x$ is conjugate to $w$,
in which case $\det(1 - q x) = \det(1 - q w)$
and the inner sum is equal to the order of
the centralizer of $w$ in $W$.  
The assertion of the theorem now follows.
\qed


\head
\leftline{\bf 4. The proof of \DihedralTheorem}
\endhead

Throughout this section we 
suppose 
$S = \left \{ r, \ s \right \}$, where
$r^2 = s^2 = (rs)^n = 1$, $n \ge 3$.  
We recall some 
facts about the
Hecke algebra $\Hecke$ of $W$
(see 
\cite{\KilmoyerSolomonRef}
or \cite{\GeckPfeifferRef, 8.3}). 
The {\Poincare} polynomial of $W$ is given by
$$
P(q) = \sum_{w \in W}q^{\ell(w)} 
= (1+q)(1+q+q^2+\cdots+q^{n-1}).
$$
For $m$ an integer, define
$\theta_m = m \pi / n$.
The assignment
$$
T_r \mapsto 
\pmatrix
-1 & 2 \cos (\theta_m) \, \sqrt{q} \cr
0 & q \cr
\endpmatrix,
\qquad
T_s \mapsto
\pmatrix
q & 0 \cr
2 \cos (\theta_m) \, \sqrt{q} & -1 \cr
\endpmatrix
$$
extends to a representation 
$\rho_m$ of $\Hecke$.
Moreover, each irreducible $2$-dimensional
representation of $\Hecke$ is equivalent to
$\rho_m$ for a unique  
$m$ satisfying $1 \le m < n/2$.
Let $\chi_m$ denote the irreducible character
of $W$ such that
$\chi_m(y)$
is obtained from $\trace \rho_m(T_y)$
by the substitution 
$\sqrt{q}\mapsto 1$.
The generic degree corresponding
to $\chi_m$ is given by
$$
d_{\chi_m}(q) =
{2 \left(1 - \cos (2\theta_m) \right) \, q \, P(q)
\over 
n \left(1 - 2 \cos( 2\theta_m) \, q + q^2\right)
} .
$$
The characteristic polynomial of
$\rho_m(T_{rs})$ is
$x^2 - 2 \cos(2 \theta_m) q x + q^2$,
so the eigenvalues of 
$\rho_m(T_{rs})$ are 
$\omega^m q$ and
$\omega^{-m} q$, where
$\omega = e^{2 m \pi i / n}$.

We require the following identity 
of formal power series.


\proclaim{\DihedralLemma}
If $a$ is an integer such that $0 \le a < n$, then 
$$
\sum\limits_{\zeta^n=1} 
{\zeta^a \over (1- \zeta q) (1 - \zeta^{-1} q) }
=
{
n ( q^{\residue{a}{n}} + q^{n-\residue{a}{n}})
\over
(1-q^2)(1-q^n) }
$$
where the sum is over all complex $n$th roots of unity
and
$\residue{a}{n}$ is the least nonnegative integer
congruent to $a$ modulo $n$.
\endproclaim

\demo{Proof of Lemma}
Define a polynomial function $f$ on $\complexes$
of degree at most $n$ by
$$
f(z) =
\sum\limits_{\zeta^n=1} 
\, \zeta^a \,
{(1-z^2)(1-z^n)
\over 
(1- \zeta z) (1 - \zeta^{-1} z) }
$$
for $z \in \complexes$, $z^n \ne 1$,
with $f$ extended continuously to $\complexes$. 
Easy calculations 
show that
$f(\zeta)=n(\zeta^{a}+\zeta^{n-a})$
whenever $\zeta^n=1$.  Also,
$f(0) = \sum_{\zeta^n=1} \zeta^a = n$ if 
$a =0$, 
and $f(0)=0$ if $0 < a < n$.
Therefore $f(z)=n(z^{a}+z^{n-a})$
for all $z \in \complexes$, 
and the assertion of the lemma follows.
\qed\enddemo 


Now, let
$w = (r s)^k$
were
$1 \le k \le n/2$.
Then 
$\chi_m(w) = \omega^{k m} + \omega^{-k m}$.
Observe
$$
g_w(q) 
= 
\sum_{\lambda} \lambda(w) d_\lambda(q)
+
\sum_{1 \le m < n/2} \chi_m(w) d_{\chi_m}(q) 
$$
where $\lambda$ ranges over the 
linear characters of $W$.

\noindent
{\it Case 1.}\ \ Suppose $n$ is odd.
In this case there are two 
linear characters 
$1$, $\varepsilon$ of $W$,
with generic degrees
$d_1(q)=1$, $d_\varepsilon(q)=q^{n}$.
Thus
$$
\eqalign{
g_w(q)
= & \ 
1 + q^{n} +
{q P(q) \over n}
\sum_{m=1}^{(n-1)/2}
(\omega^{k m} + \omega^{-k m})
{(1 - \omega^m)(1 - \omega^{-m}) 
\over (1 - \omega^m q)(1 - \omega^{-m} q)} \cr
= & \ 
1 + q^{n} +
{q P(q) \over 2 n}
\sum_{m=0}^{n-1}
(\omega^{k m} + \omega^{-k m})
{(1 - \omega^m)(1 - \omega^{-m}) 
\over (1 - \omega^m q)(1 - \omega^{-m} q)} \cr
= & \ 
1 + q^n +
{q P(q) \over 2 n}
\sum_{\zeta^n = 1}
{{2 (\zeta^{k } + \zeta^{- k }) -
     (\zeta^{k+1 } + \zeta^{-(k+1)} ) 
    - (\zeta^{k-1 } + \zeta^{-(k-1)} ) }
\over (1 - \zeta q)(1 - \zeta^{-1} q) }. \cr
}
$$
Applying {\DihedralLemma} yields  
$$
\eqalign{
g_w(q) 
= & \ 
1 + q^n +
{q P(q) \over 2n} \cdot 
{
2n \left( 2 (q^{k} + q^{n-k}) -
         (q^{k+1} + q^{n-k-1}) - (q^{k-1}+q^{n-k+1}) \right)
\over
(1-q^2)(1-q^n)
}\cr
= & \ 
1 + q^n -
{ q (1-q)^2 (q^{k-1} + q^{n-k-1}) \over (1-q)^2 }\cr
= & \ 
1 + q^n - q^k - q^{n-k} 
= (1-q^k) (1-q^{n-k}), \cr
}
$$
so the theorem is established when $n$ is odd.


\noindent
{\it Case 2.}\ \ Suppose $n$ is even.
In this case the assignments $\delta(r)=-1$, 
$\delta(s)=+1$ determine
a linear character $\delta$ of $W$. 
The four linear 
characters of $W$ are 
$1$, $\varepsilon$, $\delta$, and
$\varepsilon \delta$, with generic degrees
$$
d_1(q)=1,
\qquad
d_\varepsilon(q)=q^n,
\qquad
d_\delta(q)=
d_{\varepsilon \delta}(q)
=
{2 q (1-q^n) \over n (1-q^2)}.
$$
Hence 
$$
\eqalign{
g_w(q) 
= & \ 
1 + q^n + (-1)^k {4 q(1-q^n) \over n (1-q^2)}\cr
& \qquad\quad  +
{q P(q) \over n }
\sum_{m=1}^{(n-2)/2}
(\omega^{k m} + \omega^{- k m})
{(1 - \omega^{m})(1 - \omega^{-m}) 
\over
(1 - \omega^{m}q)(1 - \omega^{-m}q)} \cr
= & \ 
1 + q^n + (-1)^k {4 q(1-q^n) \over n (1-q^2)}
- {q P(q) \over 2 n}
(2 (-1)^k) { 4 \over (1+q)^2 } \cr
& \qquad\quad  +
{q P(q) \over 2 n }
\sum_{m=0}^{n-1}
(\omega^{k m} + \omega^{- k m})
{(1 - \omega^{m})(1 - \omega^{-m}) 
\over
(1 - \omega^{m}q)(1 - \omega^{-m}q)} \cr
= & \ 
1 + q^n + 
{q P(q) \over 2 n }
\sum_{m=0}^{n-1}
(\omega^{k m} + \omega^{- k m})
{(1 - \omega^{m})(1 - \omega^{-m}) 
\over
(1 - \omega^{m}q)(1 - \omega^{-m}q)}. \cr
}
$$
The argument given in the previous case now
applies to complete the proof of the theorem.
\qed


\head
\leftline{\bf 5. The proof of \CoxeterTheorem}
\endhead

By {\ProductProposition}, we can reduce to the case
in which $(W,S)$ is indecomposable.   
The argument presented here is organized by
cases according to the type of $(W,S)$.  
We start by 
recalling the values of the degrees 
$d_1,\dots,d_\ell$ 
of $(W,S)$
(\cite{\BourbakiRef}).

\bigskip
\tableskip
\centerline{\vbox{
\halign{\hfil # \hfil & \quad \hfil # \hfil  \cr
type of $(W,S)$  & $d_1,\dots,d_\ell$ \cr
\noalign{\vskip5pt\hrule\vskip5pt}
$A_\ell$ & 2, 3, 4, \dots, $\ell+1$ \cr
$B_\ell$, $C_\ell$ & 2, 4, 6, \dots, $2\ell$ \cr
$D_\ell$ & 2, 4, 6, \dots, $2\ell-2$, $\ell$ \cr
$E_6$ & 2, 5, 6, 8, 9, 12 \cr
$E_7$ & 2, 6, 8, 10, 12, 14, 18 \cr
$E_8$ & 2, 8, 12, 14, 18, 20, 24, 30 \cr
$F_4$ & 2, 6, 8, 12 \cr
$G_2$ & 2, 6 \cr
$H_3$ & 2, 6, 10\cr
$H_4$ & 2, 12, 20, 30\cr
$I_n$ & 2, $n$ \cr
}
\vskip5pt\hrule}}

\bigskip
\normalskip
\noindent
Here $I_n$ corresponds to the dihedral 
group of order $2n$
($n=5$ or $n \ge 7$).
The {\Poincare} polynomial of
$W$ satisfies  
$$
P(q) = \prod_{j=1}^{\ell} 
{1 - q^{d_j} \over 1-q}.
$$


\noindent{\it Type A.} \ \ 
Let $W = W(A_\ell)$.
In this case $W$ is naturally isomorphic to the
symmetric group $S_{\ell + 1}$, and the
Coxeter elements correspond to $\ell+1$ cycles.
The characteristic polynomial
of $c$ in the 
natural reflection representation is thus 
$\det(1 - q c) = (1-q^{\ell+1})/(1-q)$,
so {\TypeATheorem} gives 
$$
g_c(q) = {(1-q)^\ell P(q) \over (1-q^{\ell+1})/(1-q)}
=(1-q)(1-q^{2})\cdots(1-q^{\ell}),
$$
as required.


\noindent{\it Type B, C.} \ \ 
Suppose $W=W(B_\ell)$, where $\ell \ge 2$.   
Let $c$ be a Coxeter element of $W$:
in the natural reflection representation, $c$ 
can be taken to be the transformation  
sending $e_j$ to $e_{j+1}$ for $1 \le j < \ell$ and
$e_n$ to $-e_1$, where $e_1, \dots, e_\ell$ is 
an orthonormal 
basis for the underlying Euclidean space.
The irreducible characters 
of $W$ have the form 
$$
\chi^{\alpha,\beta}
=
\Ind_{W_j\times W_k}^W
\left(
\widetilde{\chi^\alpha}
\times
\delta
\widetilde{\chi^\beta}
\right),
$$
where $\chi^\alpha$ is the irreducible
character of
$W(A_{j-1})\cong S_j$ corresponding to
the partition $\alpha$ of $j$,
$\widetilde{\chi^\alpha}$ is the extension
of $\chi^\alpha$ to
$W_j = W(B_j)
\cong 
\left\{\pm 1\right\}^j \rtimes S_j$
with $\left\{\pm 1\right\}^j$ in the kernel,
$\widetilde{\chi^\beta}$ is defined
similarly for the partition $\beta$ of
$k = \ell - j$, and
$\delta$ is the linear character
of $W_k$ determined by
$\delta((\varepsilon_1,\dots,\varepsilon_k)\pi)
= \varepsilon_1 \cdots \varepsilon_k$
for $\varepsilon_1,\dots,\varepsilon_k
\in \left\{\pm1\right\}$, $\pi \in S_k$
(see
\cite{\GeckPfeifferRef, 5.5}).
Let $d^{\alpha,\beta}(q)$ be the
generic degree corresponding to
$\chi^{\alpha,\beta}$.
We have $\chi^{\alpha,\beta}(c)=0$ 
unless the Young diagram
of $\alpha$ is an $\ell$-hook and 
$\beta$ is empty, or vice versa.
If $\alpha=(\ell-k,1^k)$, that is, if $\alpha$
has one part $\ell-k$ and $k$ parts $1$, 
and $\beta$ is empty,
where $0 \le k \le \ell-1$,
then $\chi^{\alpha,\beta}(c)=(-1)^k$ 
and
$$
d^{\alpha,\beta}(q)
=
q^{k^2}\, {(1+q^\ell)(1+q^k) \over 2(1+q^{\ell-k})} 
\prod\limits_{j=1}^{k}{1-q^{2(\ell-j)} \over 1-q^{2j}}
$$
by the formulas of \cite{\LusztigClassicalRef} 
or \cite{\CarterFiniteGroupsRef, 13.5}.
On the other hand, if
$\alpha$ is empty and 
$\beta=(\ell-k+1,1^{k-1})$
where $1 \le k \le \ell$,
then
$\chi^{\alpha,\beta}(c)=(-1)^k$ and
$$
d^{\alpha,\beta}(q)
=
q^{k^2} \, {(1+q^\ell)(1+q^{\ell-k}) \over 2(1+q^k)}
  \prod\limits_{j=1}^{k-1} 
{1-q^{2(\ell-j)} \over 1-q^{2j}}.
$$
Notice that
$$
{1+ q^k  \over 1+ q^{\ell-k}} +
{1- q^k  \over 1- q^{\ell-k} }
= 
{2(1 - q^\ell) \over 1 - q^{2(\ell-k)}},
$$
and thus 
$$
\eqalign{
g_c(q) = &
 1 + 
\sum\limits_{k=1}^{n-1}
(-1)^k q^{k^2} 
\prod\limits_{j=1}^{k}
{1 - q^{2(\ell-j+1)} \over 1 - q^{2j}}
+ (-1)^\ell q^{\ell^2} \cr
= &
\sum\limits_{k=0}^{n}
(-1)^k q^{k^2} \, 
\qbinomial{\ell}{k}(q^2), \cr
}
$$
where 
the {\it $q$-binomial coefficient} 
$
\qbinomial{\ell}{k}(t)
$ 
is the polynomial defined by
$$
\qbinomial{\ell}{k}(t)
=
\prod\limits_{j=1}^{k}
{1 - t^{\ell-j+1} \over 1 - t^j}
=
{(1 - t^{\ell})
(1-t^{\ell-1})\cdots(1-t^{\ell-k+1})
 \over (1 - t)(1-t^2)\cdots(1-t^k)}.
$$
The 
$q$--binomial theorem or
Gauss's binomial formula 
(see
\cite{\StanleyRef, p.162} 
or
\cite{\KacCheungRef, 5.5})
gives
$$
\eqalign{
g_c(q) = &
\sum\limits_{k=0}^{\ell}
\qbinomial{\ell}{k} (q^2) 
(-q)^k \left( q^2 \right)^{\left ( {k \atop 2 } \right) }
\cr
= &
\prod\limits_{j=0}^{n-1}
(1 + (-q) (q^2)^j) \cr
= &
(1-q)(1-q^3)(1-q^5) \cdots (1-q^{2\ell-1}), \cr
}
$$
so the assertion of the theorem
holds in types $B$ and $C$.


\noindent{\it Type D.}\ \ 
Suppose $W=W(D_\ell)$, where $\ell \ge 4$.   
Let $c$ be a Coxeter element of $W$:
in the natural reflection representation, $c$ 
can be taken to be the transformation  
sending $e_j$ to $e_{j+1}$ for $1 \le j \le \ell-2$,
$e_{\ell-1}$ to $-e_1$ and
$e_\ell$ to $-e_\ell$, where $e_1, \dots, e_\ell$ is 
an orthonormal 
basis for the underlying Euclidean space.
If $\chi^{\alpha,\beta}\in \Irr(W(B_\ell))$
and $\alpha \ne \beta$, then the restriction of
$\chi^{\alpha,\beta}$ to $W$ is an irreducible character,
denoted $\chi^{\left\{\alpha,\beta\right\}}$.
Let 
$d^{\left\{\alpha,\beta\right\}}(q)$ 
be the generic degree corresponding to
$\chi^{\left\{\alpha,\beta\right\}}$.
If $\chi^{\left\{\alpha,\beta\right\}}(c) \ne 0$,
then
one of the diagrams of $\alpha$, $\beta$ must
contain an $(\ell-1)$-hook.
We can assume that the diagram
of $\alpha$ contains an $(\ell-1)$-hook
because
$\chi^{\left\{\alpha,\beta\right\}}
=\chi^{\left\{\beta,\alpha\right\}}$.
If
$\alpha=(\ell)$ and $\beta$ is empty, then
$\chi^{\left\{\alpha,\beta\right\}}=1$,
so
$\chi^{\left\{\alpha,\beta\right\}}(c)=1$
and
$d^{\left\{\alpha,\beta\right\}}(q) = 1$.
If
$\alpha=(1^\ell)$ and $\beta$ is empty,
then
$\chi^{\left\{\alpha,\beta\right\}}
= \varepsilon$,
so
$\chi^{\left\{\alpha,\beta\right\}}(c)=(-1)^\ell$
and
$d^{\left\{\alpha,\beta\right\}}(q) = 
q^{\ell^2-\ell}$.
If
$\alpha=(\ell-k-1,2,1^{k-1})$ and $\beta$
is empty, where $1\le k \le \ell-3$, then
$\chi^{\left\{\alpha,\beta\right\}}(c)=(-1)^k$
and
$$
d^{\left\{\alpha,\beta\right\}}(q) =
q^{k^2+k+1}
{(1 - q^\ell)(1 - q^{\ell-k-2}) \over
  (1 - q)(1 - q^{k+1})}
\prod\limits_{j=0}^{k} 
    {1 + q^{\ell-j-1} \over 1 + q^j}
\prod\limits_{j=1}^{k-1} 
     {1 - q^{\ell-j-1} \over 1 - q^j}
$$
by the formulas of \cite{\LusztigClassicalRef} 
or \cite{\CarterFiniteGroupsRef, 13.5}.
Finally, if
$\alpha=(\ell-k-1,1^k)$ and
$\beta=(1)$,
where $0 \le k \le \ell-2$, then
$\chi^{\left\{\alpha,\beta\right\}}(c)=(-1)^{k-1}$
and
$$
d^{\left\{\alpha,\beta\right\}}(q) =
q^{k^2+k+1}
{(1 - q^\ell)(1 + q^{k}) \over
  (1 - q)(1 + q^{\ell-k-1})}
\prod\limits_{j=0}^{k+1} 
   {1 + q^{\ell-j-1} \over 1 + q^j}
\prod\limits_{j=1}^{k} 
   {1 - q^{\ell-j-1} \over 1 - q^j}.
$$
The orthogonality relations for $W$ show 
that $\chi(c)=0$ for any $\chi \in \Irr(W)$ not
considered above.   Observe
$$
\eqalign{
{1 - q^{\ell-k-2} \over 1 - q^{k+1}} \ -& \ 
{(1 + q^{k})(1 + q^{\ell-k-2})(1 - q^{\ell-k-1}) \over
  (1 + q^{\ell-k-1})(1 + q^{k+1})(1 - q^{k})} \cr
& = \  
-2q^k 
{(1 - q)(1 - q^{\ell-1})(1 + q^{\ell-2k-2})
\over (1 - q^k)(1 - q^{2(k+1)})(1 + q^{\ell-k-1})}. \cr
}
$$
Hence for $1 \le k \le \ell-3$ we have
$$
\eqalign{
d^{\{(\ell-k-1,2,1^{k-1}),(-)\}} & (q) - 
d^{\{(\ell-k-1,1^{k}),(1)\}}(q) \cr
= & \ 
q^{k^2+k+1} 
\left({1 - q^\ell \over 1 - q}\right)
\left ( -2q^k 
{(1 - q)(1 - q^{\ell-1})(1 + q^{\ell-2k-2})
\over (1 - q^k)(1 - q^{2(k+1)})(1 + q^{\ell-k-1})} \right) \cr
& \qquad\qquad \times
\prod\limits_{j=0}^{k} 
   {1 + q^{\ell-j-1} \over 1 + q^j}
\prod\limits_{j=1}^{k-1} 
   {1 - q^{\ell-j-1} \over 1 - q^j} \cr
= & \ 
- q^{(k+1)^2}
{1 + q^{\ell-2k-2} \over 1 + q^\ell}
\prod\limits_{j=1}^{k+1} 
   {1 - q^{2(\ell-j+1)} \over 1 - q^{2j}}\cr
= & \ 
- q^{(k+1)^2}
{1 + q^{\ell-2k-2} \over 1 + q^\ell}
\qbinomial{\ell}{k+1}(q^2).
}
$$
Therefore
$$
\eqalign{
\sum\limits_{k=1}^{\ell-3}
\bigl( \bigr.
(-1)^k & d^{\{(\ell-k-1,2,1^{k-1}),(-)\}} (q)  
 +
(-1)^{k-1} d^{\{(\ell-k-1,1^{k}),(1)\}}(q) 
\bigl) \bigr. \cr
= & \ 
\sum\limits_{j=2}^{\ell-2}
(-1)^{j} q^{j^2}
{1 + q^{\ell-2 j} \over 1 + q^\ell}
\qbinomial{\ell}{j}(q^2) .\cr
}
$$
Now, as in the
argument for type $B$, we have 
$$
\eqalign{
\sum\limits_{j=0}^{\ell}
(-1)^j q^{j^2} \qbinomial{\ell}{j}(q^2)
= & \ 
(1-q) ( 1-q^3) (1-q^5) \cdots (1-q^{2\ell-1}). \cr
}
$$
Also, the $q$-binomial theorem gives
$$
\eqalign{
\sum\limits_{j=0}^{\ell}
(-1)^j q^{j^2} q^{\ell-2 j} \qbinomial{\ell}{j}(q^2)
= & \ 
q^\ell \sum\limits_{j=0}^{\ell}
\qbinomial{\ell}{j}(q^2) 
  (-q^{-1})^j 
(q^2)^{\left( {j \atop 2}\right)} 
     \cr
= & \ 
q^\ell \prod\limits_{j=0}^{\ell-1} 
\left(1 + (-q^{-1}) \cdot (q^2)^{j} \right) \cr
= & \ 
q^\ell (1-q^{-1}) (1-q) (1-q^3) \cdots (1-q^{2\ell-3}). \cr
}
$$
Hence
$$
\eqalign{
\sum\limits_{k=1}^{\ell-3} 
\bigl( \bigr. (-1)^k  &
d^{\{(\ell-k-1,2,1^{k-1}),(-)\}}(q)  +
(-1)^{k-1} d^{\{(\ell-k-1,1^{k}),(1)\}}(q) 
\bigl) \bigr. \cr
= & \ 
{1 \over 1 + q^\ell}
 (1-q)(1-q^3)(1-q^5) \cdots (1-q^{2\ell-1})  \cr
& \qquad +
{q^\ell  \over 1 + q^\ell}
(1-q^{-1})(1-q)(1-q^3) \cdots (1-q^{2\ell-3}) \cr
& \qquad -
 {1 \over 1 + q^\ell}
\sum\limits_{j \in \left \{ 0,1,\ell-1,\ell \right \}}
(-1)^j \left( q^{j^2} + q^{j^2 + \ell - 2 j} \right)
\qbinomial{\ell}{j} (q^2). \cr
}
$$
Since
$$
{(1-q^{2\ell-1}) + q^\ell (1-q^{-1}) \over 1 + q^\ell}
=
1-q^{\ell-1},
$$
it follows that 
$$
\eqalign{
\sum\limits_{k=1}^{\ell-3} 
\bigl( \bigr. (-1)^k  &
d^{\{(\ell-k-1,2,1^{k-1}),(-)\}}(q)  +
(-1)^{k-1} d^{\{(\ell-k-1,1^{k}),(1)\}}(q) 
\bigl) \bigr. \cr
= & \ 
(1-q)(1-q^3)\cdots (1-q^{2\ell-3})(1-q^{\ell-1}) - 1
+ q {(1 - q^\ell)(1 + q^{\ell-2}) \over 1 - q^2} \cr
& \qquad 
- (-1)^{\ell-1} 
   q^{(\ell-2)^2+(\ell-2)+1}
    {(1 - q^\ell)(1 + q^{\ell-2}) \over 1 - q^2} 
- (-1)^\ell q^{\ell^2-\ell}
\cr
= & \ 
(1-q)(1-q^3)\cdots (1-q^{2\ell-3})(1-q^{\ell-1}) \cr
& \qquad   
- d^{\{(\ell),(-)\}}(q) 
   + d^{\{(\ell-1),(1)\}}(q) 
   - (-1)^{\ell-1} d^{\{(1^{\ell-1}),(1)\}}(q)
   - (-1)^{\ell} d^{\{(1^\ell),(-)\}}(q). \cr
}
$$
Therefore
$$
g_c(q) =
(1-q)(1-q^3)\cdots (1-q^{2\ell-3})(1-q^{\ell-1}) ,
$$
as required.


\noindent{\it Other Types.}\ \ 
The assertion of
{\CoxeterTheorem} holds in types 
$G_2$ and $I_n$
by {\DihedralTheorem}.
The remaining cases
$E_6$, $E_7$, $E_8$, $F_4$, $H_3$ and $H_4$
have been verified by direct calculations.  The necessary
character values can be found in
\cite{\FrameOneRef},   
\cite{\FrameTwoRef},  
\cite{\KondoRef},        
\cite{\GroveRef},         
and the associated generic degrees can be found in 
\cite{\SurowskiTwoRef},                        
\cite{\BensonRef},                             
\cite{\SurowskiOneRef}, \cite{\LusztigUnipRef}, 
\cite{\LusztigClassTwoRef},                    
\cite{\AlvisLusztigRef}.                      
The details will not be shown here.
This completes the proof of the theorem.
\qed


\remark{Remark}
The list of all
polynomials $g_w(q) = \sum_{\chi} \chi(w) d_\chi(q)$
can be computed
using GAP version 3, release 4.4
(\cite{\SchonertRef})
by executing the following statements:

\tableskip
\medskip
\centerline{
\vbox{
\halign{
# \hfill \cr
RequirePackage(``chevie''); \cr
W:=CoxeterGroup({\it type},{\it rank}); \cr
chartab:=CharTable(W).irreducibles; \cr
q:=X({\it field}); \ q.name:=``q''; \cr
H:=Hecke(W,q); \cr
genericdegs:=List(SchurElements(H),
   \, x${-}{>}$PoincarePolynomial(H)/x); \cr
gpolynomials:=TransposedMat(chartab)*genericdegs; \cr
}}}

\normalskip
\noindent
Here {\it type} is the type of $(W,S)$
as a string (e.g. "E"),
{\it rank} is the rank 
$\ell = \vert S \vert$, and
{\it field} is a GAP expression for the
splitting field $K$
(e.g. Rationals in the crystallographic cases,
CyclotomicField(5) for $H_3$, $H_4$).
\endremark


\head
\leftline{\bf 6.  Concluding remarks}
\endhead

Some examples with $W$
indecomposable of small rank $\ell$ are given below.
Only elements of cuspidal classes,
that is, classes that do not meet
proper parabolic subgroups, are considered.  
In each case, the coefficients of
$g_w(q)$ are rational integers, 
$g_w(q)$ is divisible by 
$(1-q)^\ell$, and the coefficients of
the polynomial
$h_w(q)=g_w(q)/(1-q)^\ell$ are nonnegative.
The author knows of
no a priori explanation for these
phenomena.  
Also, 
$h_w(q)$ is 
palindromic by {\PalindromicProposition}, 
and so is determined by its terms  
of degree at most $(N-\ell)/2$.

In the natural reflection representation,
elements of $W(B_\ell)$ take the form of
signed permutation matrices,
that is, monomial matrices whose nonzero
entries are $\pm1$.  If 
$w\in W(B_\ell)$, the signed cycle type
of $w$ is the pair $(\mu,\nu)$, where
$\mu$ is
the partition formed by the lengths
of the cycles of $w$ with an even
number of entries $-1$ and
$\nu$ 
is the partition formed by the lengths
of the cycles with an odd
number of entries $-1$.
The conjugacy class
$C_{\mu,\nu}$ of $W(B_\ell)$ consists 
of all elements with signed cycle type
$(\mu,\nu)$.
The classes for $W(D_4)$ 
that occur in {\ExampleD} below 
 are uniquely determined
by their signed cycle types in $W(B_4)$. 
The classes of $W(F_4)$ are 
indexed by admissible graphs, as in 
\cite{\CarterClassRef}.
For type $H_3$, representatives of the
conjugacy classes of interest are given
in terms of the
generators $S=\left\{r,s,t\right\}$,
where $(rs)^3=(rt)^2=(st)^5=e$.
In each table, the first row corresponds
to the class of Coxeter elements and the
last row corresponds to the class 
containing the longest
element $w_0$ of $(W,S)$. 

\medskip\tableskip
\proclaim{\ExampleB}
$W=W(B_3)$
\endproclaim
\centerline{
\vbox{
\vskip2pt
\halign{%
 \hfil # \hfil &\quad \hfil # \hfil \cr
\text{class of $w$}  &  
$h_w(q) = g_w(q) / (1-q)^\ell $   \cr
\noalign{\vskip4pt\hrule\vskip4pt}
$C_{(-),(3)}$ & $1+2q+3q^2+3q^3+3q^4+2q^5+q^6$\cr
$C_{(-),(2,1)}$ & $1+3q+4q^2+5q^3+4q^4+3q^5+q^6$\cr
$C_{(-),(1^3)}$ & $1+2q+6q^2+6q^3+6q^4+2q^5+q^6$\cr
\noalign{\vskip 5 pt}
}
\vskip2pt\hrule\vskip2pt}}
\normalskip

\medskip
Notice that 
if $w \in W(B_3)$ 
belongs to the class
$C_{(-),(2,1)}$, then 
$g_{w}(q)$ factors as
$$
\eqalign{
g_{w}(q)
= & \ 
(1-q)^3
\left(1+q+q^2\right)
\left(1 + 2q + q^2 + 2q^3 + q^4\right), \cr
}
$$
while if $w = w_0$ belongs
to the class $C_{(-),(1^3)}$, then 
$g_{w}(q)$ factors as
$$
\eqalign{
g_{w}(q)
= & \ 
(1-q)^3
\left(1+q+q^2\right)
\left(1 + q + 4q^2 + q^3 + q^4\right). \cr
}
$$
Hence $g_w(q)$ is not in general
a product of cyclotomic polynomials.

\medskip\tableskip
\proclaim{\ExampleD}
$W=W(D_4)$
\endproclaim
\centerline{
\vbox{
\vskip2pt
\halign{%
 \hfil # \hfil &\quad \hfil # \hfil \cr
\text{class of $w$}  &  
$h_w(q) = g_w(q)/(1-q)^\ell$   \cr
\noalign{\vskip4pt\hrule\vskip4pt}
$C_{(-),(3,1)}$ & 
$1+3q+6q^2+8q^3+9q^4+O(q^5)$\cr
$C_{(-),(2,2)}$ & 
$1+4q+7q^2+10q^3+10q^4+O(q^5)$\cr
$C_{(-),(1^4)}$ & 
$1+3q^2+8q^3+12q^4+O(q^5)$\cr
\noalign{\vskip 5 pt}
}
\vskip2pt\hrule\vskip2pt}}
\normalskip

\medskip\tableskip
\proclaim{\ExampleF}
$W=W(F_4)$
\endproclaim
\centerline{
\vbox{
\vskip2pt
\halign{%
 \hfil # \hfil &\  \hfil # \hfil \cr
\text{class of $w$}  &  
$h_w(q) = g_w(q)/(1-q)^\ell$   \cr
\noalign{\vskip4pt\hrule\vskip4pt}
$F_4$ &
$1+3q+6q^2+10q^3+15q^4+20q^5+25q^6+
29q^7+32q^8+34q^9+35q^{10}+O(q^{11})$\cr
$F_4(a_1)$ &
$1+4q+10q^2+16q^3+24q^4+32q^5+40q^6
+44q^7+50q^8+52q^9+54q^{10}+O(q^{11})$\cr
$B_4$ &
$1+4q+9q^2+16q^3+24q^4+33q^5+41q^6
+48q^7+53q^8+57q^9+58q^{10}+O(q^{11})$\cr
$D_4$ & 
$1+4q+10q^2+19q^3+30q^4+41q^5+52q^6
+62q^7+71q^8+76q^9+78q^{10}+O(q^{11})$\cr
$D_4(a_1)$ &
$1+6q+15q^2+28q^3+43q^4+60q^5+75q^6
+86q^7+95q^8+104q^9+108q^{10}+O(q^{11})$\cr
$C_3+A_1$ &
$1+4q+10q^2+19q^3+30q^4+41q^5+52q^6
+62q^7+71q^8+76q^9+78q^{10}+O(q^{11})$\cr
$A_3+\widetilde{A}_1$ &
$1+3q+7q^2+14q^3+24q^4+36q^5+
49q^6+63q^7+75q^8+82q^9+84q^{10}+O(q^{11})$\cr
$A_2+\widetilde{A}_2$ &
$1+2q+2q^2+4q^3+9q^4+18q^5+32q^6+46q^7+54q^8
+58q^9+60q^{10}+O(q^{11})$\cr
$4A_1$ & 
$1+4q+19q^2+40q^3+66q^4+80q^5+97q^6
+116q^7+155q^8+184q^9+204q^{10}+O(q^{11})$\cr
\noalign{\vskip 5 pt}
}
\vskip2pt\hrule\vskip2pt}}
\normalskip

\medskip\tableskip
\proclaim{\ExampleH}
$W=W(H_3)$
\endproclaim
\centerline{
\vbox{
\vskip2pt
\halign{%
 \hfil # \hfil &\quad \hfil # \hfil \cr
\text{$w$}  &  
$h_w(q) = g_w(q) / (1-q)^\ell $   \cr
\noalign{\vskip4pt\hrule\vskip4pt}
$rst$ & $1+2q+3q^2+4q^3+5q^4+5q^5+5q^6+O(q^7)$\cr
$rstst$ & $1+3q+5q^2+7q^3+8q^4+9q^5+9q^6+O(q^7)$\cr
$rstsrstst$ & $1+3q+6q^2+9q^3+12q^4+14q^5+15q^6+O(q^7)$\cr
$w_0$ & $1+2q^2+4q^3+11q^4+12q^5+15q^6+O(q^7)$\cr
\noalign{\vskip 5 pt}
}
\vskip2pt\hrule\vskip2pt}}
\normalskip

Based on these and
other examples with small rank, 
we state
the following conjecture.  

\proclaim{\Conjecture}
If $m$ is the multiplicity of $1$
as an eigenvalue of $w \in W$ in the natural 
reflection representation of $W$,
then
$$
g_w(q) 
=
(1-q)^{\ell - m}
h_w(q)
$$
where $h_w(q)$ 
is a palindromic polynomial 
whose coefficients are nonnegative rational integers.
\endproclaim

By {\ProductProposition} and
{\RecursionProposition},  the conjecture holds
in general if it holds whenever 
$(W,S)$ is indecomposable and
$w$ belongs to a cuspidal class.  
(For $J \subseteq S$, it is known that 
$P(q)/P_{W_J}(q) =
\sum_{x \in D_J} q^{\ell(x)}$, where
$D_J$ is the set of minimal length left coset
representatives of $W_J$ in $W$.
In particular, $P(q)/P_{W_J}(q)$ has
nonnegative coefficients.)
The conjecture holds in type $A$
by {\TypeATheorem}, and holds
in types $G_2$ and $I_n$ by
{\DihedralTheorem}.
The conjecture has been verified by direct
calculations for cuspidal classes 
in types $E_6$, $E_7$, 
$E_8$, $F_4$, $H_3$ and $H_4$.
Hence to verify the conjecture, 
it is enough to establish the cases
$W=W(B_\ell)$ and
$W=W(D_\ell)$ when $w$ belongs to
a cuspidal class.


\enddocument
\bye